\documentclass[a4paper,12pt]{amsart}
\usepackage{amsfonts}
\usepackage{amsthm}
\usepackage{amssymb}
\usepackage{amsmath}
\usepackage{mathabx}
\usepackage{tikz-cd}

\def\bbC{{\mathbb C}}
\def\bbD{{\mathbb D}}

\def\bbN{{\mathbb N}}

\def\cD{{\mathcal D}}

\def\cH{{\mathcal H}}

\def\cR{{\mathcal R}}

\def\om{{\omega}}
\def\whd{{\widehat{\cD}}}
\def\wcd{{\widecheck{\cD}}}
\def\whom{{\widehat{\omega}}}

\def\gtr{{\gtrsim}}

\def\awp{{A_\omega^p}}
\def\avp{{A_v^p}}

\usepackage[colorlinks=blue]{hyperref}
\usepackage{graphics}
\usepackage{hyperref}
\hypersetup{
	colorlinks=true,
	linkcolor=red,
	citecolor=cyan,
}
\newtheorem{remark}{Remark}

\newtheorem{theorem}{Theorem}[section]

\newtheorem{lemm}[theorem]{Lemma}
\newtheorem{lemma}{Lemma}

\title[Two weights inequality for Hankel operators]{Two weights inequality for Hankel operators on weighted Bergman spaces induced by radial weights}
\author{Mingjin Li, Jianren Long$^*$, Pengcheng Wu}
\date{}
\address{Mingjin Li \newline
	School of Mathematical Sciences, Guizhou  Normal University, Guiyang, 550025,  China.}
\email{limingjin2022@163.com }
\address{Jianren Long \newline
	School of Mathematical Sciences, Guizhou  Normal University, Guiyang, 550025,  China.}
\email{longjianren2004@163.com }
\address{Pengcheng Wu \newline
	School of Mathematical Sciences, Guizhou  Normal University, Guiyang, 550025,  China.}
\email{wupc@gznu.edu.cn}
\thanks{*Corresponding author. }
\begin{document}
	\maketitle
	\begin{abstract}
		The two weights inequality for Hankel operators 
		\[\|H_f^\om (\cdot)\|_{L_\eta^q}\leq C \|\cdot\|_{A_v^p}, \] induced by some radial weights under the regular assumptions  is considered, the boundedness and compactness of Hankel operators \(H^\omega_f\) is characterized for \(1<p,q<\infty\) and \(\om\neq v\neq \eta\in \cR\).
	\end{abstract} 
	\section{introduction}
	A radial weight \(\om\) belongs to \(\whd\) if $\whom(z)=\int_r^1\om(s)ds$ satisfies the doubling conditions, i.e. there exists $C=C(\om)\geq 1$ such that $$\whom(r)\leq C\whom(\frac{1+r}{2}), \quad 0\leq r<1.  $$
	Furthermore, if there exist $K=K(\omega)>1$ and $C=C(\omega)>1$ such that \[\whom(r)\ge C\whom\left(1-\frac{1-r}{K}\right),\quad 0\le r<1,\]
	then we write $\omega\in\wcd$. The definitions of both $\whd$ and $\wcd$ have their obvious geometric interpretations. The intersection $\whd\cap\wcd$ is denoted by $\cD$. For \(\alpha>-1\),  \(\om(z)=(1+\alpha)(1-|z|^2)^\alpha\) is a classical \(\cD\) function,  so the class \(\cD\) naturally arises in operator theory on spaces of analytic functions,  For \(\om\in\cD\) belong to \(\cR\) if $$\whom(r)\asymp\whom(r)(1-r), 0\leq r<1.$$
	The class \(\cD\) and \(\cR\) is the largest possible class of functions for which most of the results concerning the standard Bergman space \(A_\alpha^p\) holds \cite[Page 2]{Duan2022}. For more properties of radial weights, it can be found in \cite{Pelaez2014,Pelaez2016summer,Pelaez2018} and therein references. 
	Let $\bbD=\{z\in\bbC:|z|<1\}$ be the unit disc and  $\cH(\bbD)$ denotes the set of the all holomorphic functions on $\bbD$. 
	
	Given some nonnegative functions \(\om\) in \(\bbD\), for \(0<p<\infty\), the weighted Lebesgue spaces \(L_\om^p\) consists of all measurable functions with 
	\[\|f\|_{L_\om^p}=\left(\int_\bbD|f(z)|^p\om(z)dA(z)\right)^\frac{1}{p}<\infty,\] where \(dA(z)=dxdy=\frac{r}{\pi}drd\theta\) is the normalized area measure on \(\bbD\).
	Weighted Bergman spaces \(\awp\) is the closed spaces of \(L_\om^p\), consists of all \(f\in\cH(\bbD)\) with norm \(\|f\|_{L_\om^p}\).  If \(\om(z)=(1-|z|^2)^\alpha,\alpha>-1\), then \(A_\om^p\) is the standard radial weight Bergman spaces, we write \(A_\alpha^p\) for simple. If \(\alpha=0\) or  \(\om(z)=1\), it is the standard Bergman space \(A^p\).
	
	A series of works by Pel\'aez et al. extended some results from harmonic analysis and operator theory to a larger class of function spaces, such as Carleson measures\cite{Pelaez2015}, Hilbert operators\cite{Pelaez2013}, Bergman projections \cite{Pelaezprojection2016, Pelaez2021} and so on,  which was extended to \(A_\om^p\) instead of \(A_\alpha^p\),  in \cite{Pelaezprojection2016} and \cite{Pelaez2021} so called two weights inequality for Bergman projection 
	\[\|P_\om (\cdot)\|_\eta\leq C\|\cdot\|_v\] was characterized for \(\eta=v\). The study of two weights inequality for classical operators have attracted a considerable amount of attention and closely to other interesting questions, see \cite{Aleman2017,Hytonen2012,Korhonen2021,Lacey2014,Pelaezprojection2016,Pelaez2021,Zhang2025} and therein references for details.
	
	For \(\om\) is radial, point evaluation functional is continuous on Hilbert spaces \(A_\om^2\), hence by Riesz representation theorem, there exists a unique function \(B_z^\om(\zeta)=B^\om(z,\zeta)=\overline{B_\zeta^\om(z)}\) such that 
	\[\delta_z(f)=f(z)=\langle f, B_z^\om\rangle_{A_\om^2}.\] The function \(B_z^\om\) is called Bergman kernel such that the orthogonal Bergman projection \(P_\om\) from \(L_\om^2\) to \(A_\om^2\) is defined as 
	\[P_\om(f)(z)=\int_\bbD f(\zeta)B^\om(z,\zeta)\om(\zeta)dA(\zeta).\]
	
	Suppose \(f\in L_\om^1\), Hankel operator is defined as \[H^\om_f(g)=fg-P_\om(fg)=(Id-P_\om)(fg),\] where \(Id\) is the identity operator.
	For \(\om\) is radial, the dilated function \(f_r(z)=f(rz)\) approximate the function \(f\in\awp\), this implies that polynomials are dense in \(A_\om^p\), so \(H^\infty(\bbD)\) is dense in \(A_\om^p\), where \(H^\infty(\bbD)\) denotes all bounded analytic functions in \(\bbD\). Hence, \(H_f^\om\) is well defined in the sense of density. In what follows, if \(\om\equiv 1\) or \(\om=(1+\alpha)(1-|z|^2)^\alpha\), then we write \(H_f^\om\) as \(H_f\) for simple.

	The study of Hankel operators on Bergman spaces started in \cite{Axler1986} which say that \(H_{\overline{f}}\) is bounded or compact from \(A^2\) to \(L^2\) if and only if \(f\) is in Bloch or little Bloch spaces. Bekoll\'e et al. \cite{Bekolle1990} researched the same problem on bounded symmetric domains but with the restriction that \(H_f\) and \(H_{\overline{f}}\) are simultaneously bounded. A new method was used by  
	Luecking \cite{Luecking1992}, which characterized those function \(f\in L^p\) such that \(H_f:A^p\to L^p\) is bounded by bounded distance to analytic (BDA) functions and vanshing distance to analytic (VDA) functions, At the same time, Luecking also mentioned that the approach can only work with \(-1<\alpha<\frac{1}{p-1}\) for standard weighted Bergman spaces \(A_\alpha^p\). Hence,  he raised the question:  Whether similar characterizations exist for weighted Bergman spaces?  Pau\cite{Pau2013} characterized Hankel operators on \(A_\alpha^p\) for all \(p>1\) and  \(\alpha>-1\). For other  study of Hankel operators has attracted widespread attention from numerous scholars, see  \cite{Arazy1988,Hu2018,Hu2023,Jason1988,Stroethoff1990,Pau2016,Pelaez2020,Duan2023} and therein references for details. 
	
	However,  two weights inequality for Hankel operators
	\begin{equation}\label{twoweights}
		\|H_f^\om (\cdot)\|_{L_\eta^q}\leq C \|\cdot\|_{A_v^p}, 
	\end{equation}
	there are few literature with \(p\neq q\) even if \(\om=v=\eta\equiv 1\). Recently, for Hankel operators \(H_f^\om:A_\om^p\to L_\om^q\), Hu et al. \cite{Hu2019} characterized the boundedness and compactness of Hankel operators for the case of \(\om\equiv v\equiv\eta\) by using BDA and VDA functions for all \(1<p,q<\infty\).  Motivated by this, two weights inequality for Hankel operators \eqref{twoweights} will be characterized under some regular assumptions in this paper.
	Before give our main results, we need more notations. Let \[\beta(z,w)=\frac{1}{2}\log\frac{1+|\varphi_z(w)|}{1-|\varphi_z(w)|} \] be the Bergman metric, where \(\varphi_z(w)=\frac{z-w}{1-\overline{z}w}\), and \(D(z,r)=\{w:\beta(z,w)<r\}\) be the Bergman disc with center \(z\) and radius \(r\). A sequence \(\{a_k\}\) in \(\bbD\) is called an \(r\)-lattice (\(r > 0\)) in the Bergman metric if the following conditions are satisfied: 
	\begin{enumerate}
		\item [(i)] \(\bbD=\bigcup\limits_{k=1}^\infty D(a_k, r)\),
		\item [(ii)] \(\beta(a_i, a_j) \geq \frac{r}{2}\) for all \(i\) and \(j\) with \(i \neq j\).
	\end{enumerate}
	By \cite[Chapter 4]{Zhu2007}, for any \(r>0\), there exists an \(r-\)lattice in the Bergman metric, and there exists a positive constant \(N\) such that every point \(z\in\bbD\) belong to at most \(N\) disks \(D(a_k,r)\). 
	
	This paper is organized as follows. Some preliminary results are given in Section 2, The characterization for two weights inequality of  Hankel operators induced by radial weights is given in Section \(3\) for the case of \(1<p\leq q<\infty\), and for the case of \(1<q<p<\infty\) is shown in Section \(4\).
	Throughout this paper, \(C=C(\cdot)\) denotes the positive constant whose values may different, we write \(A\lesssim B\) or \(A\gtrsim B\) if there exist some \(C\) such that \(A\leq CB\) or \(A\geq CB\), we write \(A\asymp B\) if \(A\lesssim B\) and \(A\gtrsim B\).
	\section{Preliminary results}
	The properties of the regular radial weights will be used which are summary results from \cite{Pelaez2014,Pelaez2016summer,Pelaez2018}.
	\begin{lemma}\label{weights}
		Let \(\omega\) be a radial weight. 
		\begin{itemize}
			\item[(i)] \(\om\in\cD\) if and only if there exist \(\beta = \beta(\omega) > 0\) such that
			\[
			\widehat{\omega}(r) \asymp \left(\frac{1 - r}{1 - t}\right)^{\beta} \widehat{\omega}(t), \quad 0 \leq r \leq t < 1;
			\]
			\item [(ii)] If \(\om\in \cR\), then \(\om(D(z,r))\asymp \widehat{\omega}(z)(1 - |z|)\), \(|z| \to 1^{-}\); For \(\zeta\in D(z,r)\), \(\whom(z)(1-|z|)\asymp \whom(\zeta)(1-|\zeta|)\).
			\item [(iii)] Let \( \omega \in \widehat{\cD} \). Then there exists \( r = r(\omega) \in (0,1) \) such that 
			\[
			|B_a^\omega(z)| \asymp B_a^\omega(a)
			\] 
			for all \( a \in \mathbb{D} \) and \( z \in \Delta(a,r) \) where \(\Delta(a,r)=\{z\in\bbD:\rho(z,a)<r\}\) and \(\rho(z,a)=\left|\frac{z-a}{1-\overline{a}z}\right|\).
		\end{itemize}
	\end{lemma}
	\begin{proof}
		The proof of (i) can be found in \cite{Pelaez2016summer} and (ii) can be found in \cite[chapter 1]{Pelaez2014}, (iii) can be found in \cite[Lemma 8]{Pelaez2018}.
	\end{proof}
	The following lemma plays an important role in the estimates of Bergman kernel. 
	\begin{lemma}\cite[Theorem 1]{Pelaezprojection2016}\label{normestimates}
		Let \(\omega, v \in \widehat{D}\), \(0 < p < \infty\), and \(n \in \mathbb{N}\). Then
		
		\[
		\|(B_z^{\omega})^{(n)}\|_{A_v^p}^p \asymp \int_0^{|z|} \frac{\widehat{v}(t)}{\widehat{\omega}(t)^p (1-t)^{p(n+1)}} \, dt, \quad |z| \to 1^-.
		\]
		In particular, if \(1 < p < \infty\), \(\omega \in \mathcal{R}\), and \(r \in (0, 1)\), then
		\[
		\|B_z^{\omega}\|_{A_\omega^p}^p  \asymp \frac{1}{\omega(\Delta(z, r))^{p-1}}, \quad z \in \mathbb{D}.
		\]
	\end{lemma}
	\begin{remark}\label{rem1}
		By Lemmas \ref{weights} and \ref{normestimates}, 
		for \(\om, v\in\cR\), 
		\begin{align*}
			\|B_z^\om\|_{\avp}&\asymp\left(\int_0^{|z|} \frac{\widehat{v}(t)}{\widehat{\omega}(t)^p (1-t)}dt\right)^\frac{1}{p}\\
			&\asymp\left(\frac{\frac{\whom(z)}{(1-|z|)^\beta}}{\frac{\whom(z)^p}{(1-|z|)^{p\beta}}}\int_0^{|z|}\frac{dt}{(1-t)^{p\beta+p-\beta}}\right)^\frac{1}{p}\\
			&\asymp \frac{(\widehat{v}(z)(1-|z|))^\frac{1}{p}}{\whom(z)(1-|z|)}\\
			&\asymp\frac{v(D(z,r))^\frac{1}{p}}{\om(D(z,r))}.
		\end{align*}
	\end{remark}

	Let \(\mu\) be a positive Borel measure on \(\bbD\) and \(0<q<\infty\). \(\mu\) is called a \(q\)-Carleson measure for \(A_\om^p\) if formal identity operator \(Id: A_\om^p\hookrightarrow L^q(\mu)\) is bounded, and \(\mu\) is called a vanishing \(q\)-Carleson  measure if \(Id: A_\om^p\hookrightarrow L^q(\mu)\) is compact. Pel\'aez et al. obtained the following two results.
	\begin{lemm}\cite[Theorem 1]{Pelaez2015}\label{CMpq}
		Let \( 1 < p\leq q < \infty \) and \(\omega \in \cR\), and let \(\mu\) be a positive Borel measure on \(\mathbb{D}\).
		\begin{itemize}
			\item [(i)] $\mu$ is a $q$-Carleson measure for $A^{p}_{\omega}$ if and only if
			\[
			\sup_{z \in D} \frac{\mu(D(z,r))}{\omega(D(z,r))^{\frac{q}{p}}} < \infty
			\]
			for some (or any) $r > 0$. Moreover,
			\[
			\|Id\|^{q}_{A^{p}_{\omega} \to L^{q}(d\mu)} \asymp \sup_{z \in D} \frac{\mu(D(z,r))}{\omega(D(z,r))^{\frac{q}{p}}}.
			\]
			
			\item [(ii)] $\mu$ is a vanishing $q$-Carleson measure for $A^{p}_{\omega}$ if and only if
			\[
			\lim_{|z| \to 1} \frac{\mu(D(z,r))}{\omega(D(z,r))^{\frac{q}{p}}} = 0
			\]
			for some (or any) $r > 0$.
		\end{itemize}
	\end{lemm}

	\begin{lemm}\cite[Theorem 1]{Pelaez2015}\label{CMqp}
		Let \( 1 < q < p < \infty \) and \(\omega \in \mathcal{R}\), and let \(\mu\) be a positive Borel measure on \(\mathbb{D}\). Then the following statements are equivalent.
		\begin{enumerate}
			\item [(i)] \(\mu\) is a \(q\)-Carleson measure for \(A^p_{\omega}\).
			
			\item [(ii)] \(\mu\) is a vanishing \(q\)-Carleson measure for \(A^p_{\omega}\).
			
			\item [(iii)] For some (or any) \(r > 0\), \[\frac{\mu(D(\cdot, r))}{\omega(D(\cdot, r))} \in L^{\frac{p}{p-q}}_{\omega}.\]
		\end{enumerate}
		Furthermore,
		\[
		\|Id\|^q_{A^p_{\omega} \to L^q(d\mu)} \asymp \left\|\frac{\mu(D(\cdot, r))}{\omega(D(\cdot, r))}\right\|_{L^{\frac{p}{p-q}}_{\omega}}.
		\]
	\end{lemm}
	\begin{remark}
		The above two results in \cite{Pelaez2015} is obtained by using the pseudohyperbolic distance \(\rho(z,w)\) and the assumption \(\om\in\widehat{\cD}\).
	\end{remark}
	For an \(r\)-lattice \(\{a_j\}_{j=1}^\infty\), let \(\phi_j\) be some partition of unity subordinate to \(D(a_j,r)\), i.e. 
	\[\phi_j\in C^\infty(\bbD),\,\, \text{Supp}\,\phi_j\subset D(a_j,r),\,\, \phi_j\geq 0,\,\, \sum_{j=1}^\infty \phi_j=1.\]

	For a function \(f\in C^k(\bbD)\) and \(g\in\cH(\bbD)\), \(\overline{\partial}(fg)=g\overline{\partial}f\), so the boundedness of Hankel operators equivalent to the existences of the (weak) solutions of a \(\overline{\partial}\)-equation with optimal regular estimates. 
	\begin{lemm}\cite[Lemma 4.1]{Hu2019}\label{dbar}
		Suppose \(\omega \in \mathcal{R}\), \(\phi_j\) be some partition of unity subordinate to \(D{(z_j,r)}\), \(f \in L_{\omega}^{1}\) and \(\overline{\partial} f \in L^p_{(1-|\cdot|)\omega}\) for some \(p > 1\), where the derivative is in the sense of distribution. Then for \(g \in H^{\infty}\) there holds
		\[
		H_f(g) = u - P(u),
		\]
		where
		\[
		u(z) = \sum_{j=1}^{\infty} B_{z_j}(z) \int_{\bbD} \frac{\phi_j(\xi)}{(\xi - z)B_{z_j}(\xi)} g(\xi) \overline{\partial} f(\xi) dA(\xi).
		\]
	\end{lemm}
	We shall need the following atomic type results for the case of \(1<q<p<\infty\).
	\begin{lemm}\cite[Theorem 3]{martin2021}\label{duality}
		Let \(1 < p < \infty\), \(\frac{1}{p}+\frac{1}{p'}=1\) and \(\omega \in \widehat{\mathcal{D}}\), and let \(\nu\) be a radial weight. Then \(P_\omega : L^p_\nu \to L^p_\nu\) is bounded if and only if \(\left( A^{p'}_\sigma \right)^* \simeq A^p_\nu\) via the \(A^2_\omega\)-pairing with equivalence of norms, where \(\sigma=\left(\frac{\om}{v^\frac{1}{p}}\right)^{p'}\).
	\end{lemm}
	By \cite[Theorem 13]{Pelaez2021}, for \(p>1\), \(\om,v\in{\cD}\),  \(P_\om: L_v^p\to L_v^p\) is bounded if and only if 
	\[A_p(\om,v)=\sup_{0\leq r<1}\frac{\widehat{v}(r)^\frac{1}{p}\widehat{\sigma}(r)^\frac{1}{p'}}{\widehat{\om(r)}}<\infty.\]  It follows from this fact and Lemma \ref{duality} that the following results is obtained.
	\begin{lemm}\label{atmoic}
		Suppose that \(p>1\), \(\sigma=\left(\frac{\om}{v^\frac{1}{p}}\right)^{p'}\), \(\frac{1}{p}+\frac{1}{p'}=1\), \(\sigma=\left(\frac{\om}{v^\frac{1}{p}}\right)^{p'}\), \(\om\in\cR\), \(v\in\cR\), \(\eta\in\cR\), \(\{\lambda_j\}_{j=1}^\infty\in l^p\). Set \[F(z)=\sum_j\lambda_j\frac{B_{z_j}^\om(z)}{\|B_{z_j}^\om\|_{A_v^p}}.\] Then \(\|F(z)\|_{A_v^p}\lesssim 1.\)
	\end{lemm} 
	
	\begin{proof}
		Suppose that  \(\{z_j\}\) is an \(r\)-lattice in \(\bbD\) , then Lemma \ref{weights}, Remark \ref{rem1} and H\"older's inequality yields
			\begin{align*}
				& \left|  \langle g,F\rangle_{A_\om^2}\right|\\
				&=\left|\sum_j\overline{\lambda_j}\frac{g(z_j)}{\|B_{z_j}\|_{A_v^p}}\right|\\
				&\leq\|\lambda_j\|_{l^p}\left(\sum_j\frac{|g(z_j)|^{p'}}{\|B_{z_j}^\om\|_{A_v^p}^{p'}}\right)^\frac{1}{p'}\\
				&\asymp\|\lambda_j\|_{l^p}\left(\sum_j\left(\frac{\om(D(z_j,r))}{v(D(z_j,r))^\frac{1}{p}}\right)^{p'}|g(z_j)|^{p'}\right)^\frac{1}{p'}\\
				&\lesssim \|\lambda_j\|_{l^p}\left(\sum_j\left(\frac{\om(D(z_j,r))}{v(D(z_j,r))^\frac{1}{p}}\right)^{p'}\frac{1}{\left(\frac{\om(D(z_j,r))}{v(D(z_j,r))^\frac{1}{p}}\right)^{p'}}\int_{D(z_j,r)}|g(z)|^{p'}\left(\frac{\om}{v^\frac{1}{p}}\right)^{p'}dA(z)\right)^\frac{1}{p'}\\
				&\lesssim\|\lambda_j\|_{l^p}\left(\int_\bbD|g(z)|^{p'}\left(\frac{\om}{v^\frac{1}{p}}\right)^{p'}dA(z)\right)^\frac{1}{p'}. 
			\end{align*}
			Since \(A_p(\om,v)<\infty\), then for \(g\in A_\sigma^{p'}\), Lemma \ref{duality} and an argument of duality yields 
			\[\|F\|_{A_v^p}=\sup_{\|g\|_{A_\sigma^{p'}}\leq 1} \langle g,F \rangle_{A_\om^2}\lesssim 1.\]
				The proof is completed.
		\end{proof}
		
		\section{The case of \(1<p\leq q<\infty\)}
		The average function \(M_r(f)(z)=\frac{1}{|D(z,r)|}\int_{D(z,r)}fdA\) and the BDA function is defined as
		\[G_{q,r}(f)(z)=\inf\left\{\left(\frac{1}{D(z,r)}\int_{D(z,r)}|f-h|^qdA(z)\right)^\frac{1}{q}: h\in\cH(D(z,r))\right\}.\] For \(v\in\cR\),
		\[G_{q,r}(f)(z)\asymp\inf\left\{\left(\frac{1}{v(D(z,r))}\int_{D(z,r)}|f-h|^qvdA(z)\right)^\frac{1}{q}: h\in\cH(D(z,r))\right\}.\]
		For \(\cdot\in\cR\), we write  \[[\cdot]=\widehat{\cdot}(z)(1-|z|)\asymp\cdot(D(z,r)).\]

		Now we give  main results of this section.
		\begin{theorem}\label{mainth1}
			Let \(\om,v,\eta\in\cR\), \(1<p\leq q<\infty\), \(A_p(\om,\eta)<\infty\), \(f\in L_v^1\). Then the following statements are equivalent. 
			\begin{enumerate}
				\item [(i)] \(H_f^\om: A_v^p\to L_\eta^q\) is bounded;
				\item [(ii)] For some (or any) \(r>0\), \([\eta]^\frac{1}{q}[v]^{-\frac{1}{p}}G_{q,r}(f)(z)<\infty\);
				\item [(iii)] \(f\) admits a decomposition \(f=f_1+f_2\), where \(f_1\in C^1(\bbD)\)  satisfies 
				\[[\eta]^\frac{1}{q}[v]^{-\frac{1}{p}}(1-|z|)|\overline{\partial}f_1|<\infty,\]
				and \(f_2\) satisfies 
				\[[\eta]^\frac{1}{q}[v]^{-\frac{1}{p}}M_r(|f_2|^q)^\frac{1}{q}<\infty.\]
			\end{enumerate}
		\end{theorem}
		\begin{remark}\label{remark3}
			If \(\om=v=\eta\), then Theorem \ref{mainth1} implies \cite[Theorem 4.2]{Hu2019}.
		\end{remark}
		
		\begin{proof}
			We prove the theorem in its natural order. 
			
			\((\text{i})\implies (\text{ii})\).
			Suppose that \(H^\om_f: A_v^p\to L_\eta^q\) is bounded. Let \(b_z^\om=\frac{B_z^\om}{\|B_z^\om\|_{A_\om^2}}\), then by Lemma \ref{weights}, for \(w\in D(z,r)\), 
			\[|B_z(w)|\asymp |B_z(z)|\asymp\frac{1}{\whom(z)(1-|z|)}.\]
			On the one hand, 
			\[|b_z^\om|\gtr\frac{1}{(\whom(z)(1-|z|))^\frac{1}{2}}, \quad w\in D(z,r).\]
			\begin{align}\label{eq2}
				\begin{aligned}
					\|H^\om_f (b_z^\om)\|_{L_\eta^q}^q&=\int_\bbD\left| fb_z^\om-P_\om(fb_z^\om)\right|^q\eta dA\\
					&\geq\int_{D(z,r)}|b_z^\om|^q\left| f-\frac{P_\om(fb_z^\om)}{b_z^\om}\right|^q\eta dA\\
					&\gtrsim \frac{1}{\whom(z)(1-|z|)^\frac{q}{2}}{\widehat{\eta}(z)}{(1-|z|)}G_{q,r}^q(f)(z)\\
					&=[\om]^{-\frac{q}{2}}[\eta]G_{q,r}^q(f)(z).
				\end{aligned}
			\end{align}
			On the other hand, 
			\begin{align}
				\begin{aligned}\label{eq3}
					\|H_f^\om(b^\om_z)\|^q_{L_\eta^q}\leq \|H_f^\om\|^q\|b_z^\om\|_{A_v^p}^q\lesssim[v]^\frac{q}{p}[\om]^{-\frac{q}{2}}.
				\end{aligned}
			\end{align}
			The assumption with \eqref{eq2} and \eqref{eq3} yields 
			\begin{equation}\label{eqHf}
				[\eta]^\frac{1}{q}[v]^{-\frac{1}{p}}G_{q,r}(f)(z)\lesssim\|H_f^\om\|<\infty. 
			\end{equation}
			
			\(\text{(ii)}\implies \text{(iii)}\).
			According to \cite[Chapter 8]{Zhu2007}, \(G_{q,r}(f)(z)\) is essentially independent of the choice of \(r\).  For an \(\frac{r}{2}\)-lattice \(\{a_k\}\), there exists a smooth partition of unity subordinate to \(D(a_j,r)\) with \((|1-|a_j|)|\overline{\partial}\phi_j|\leq C\). By the definition of \(G_{q,r}(f)(a_j)\), there exists \(h_{j,k}\in\cH(D(a_j,r))\) such that 
			\begin{align}\label{eqhjk}
				\begin{aligned}
					G_{q,r}(f)(a_j)&\leq \left(\frac{1}{v(D(z,r))}\int_{D(z,r)}|f-h|^qvdA\right)^\frac{1}{q}\\
					&\leq G_{q,r}(f)(a_j)+\frac{1}{k},\,\,k\in\bbN^+.
				\end{aligned}
			\end{align}
			this implies that
			\begin{align*}
				|h_{j,k}|^q&\lesssim  \frac{1}{v(D(a_j,r))}\int_{D(a_j,r)}|h_{j,k}|vdA\\
				&\leq\left(\frac{1}{v(D(a_j,r))}\int_{D(a_j,r)}|h_{j,k}|^qvdA\right)^\frac{1}{q}\\
				&\leq\left(\frac{1}{v(D(a_j,r))}\int_{D(a_j,r)}|f-h_{j,k}|^qvdA\right)^\frac{1}{q}\\
				&+\left(\frac{1}{v(D(a_j,r))}\int_{D(a_j,r)}|f|^qvdA\right)^\frac{1}{q}\\
				&\leq G_{q,r}(f)(a_j)+\frac{1}{m}+\left(\frac{1}{v(D(a_j,r))}\int_{D(a_j,r)}|f|^qvdA\right)^\frac{1}{q},\,\,m\in\bbN^+.
			\end{align*}
			This implies that \(\{h_{j,k}\}_{k=1}^\infty\) is normal family functions, hence there exists a function \(h_j\in\cH(D(z,r))\) such that \(h_{j,k}\rightrightarrows h_j\) on the compact subset of \(D(a_j,r)\). Applying Fatou's lemma to \eqref{eqhjk}, we get 
			\begin{equation}
				G_{q,r}(f)(a_j)=\left(\frac{1}{v(D(a_j,r))}\int_{D(a_j,r)}|f-h_j|^qvdA\right)^\frac{1}{q}.
			\end{equation}
			For \(z\in D(a_i,r)\bigcap D(a_j,r)\), 
			\begin{align}\label{eqhihj}
				\begin{aligned}
					|h_i-h_j|&\lesssim\left(\frac{1}{v(D(z,r))}\int_{D(z,r)}|h_i-h_j|^qvdA\right)^\frac{1}{q}\\
					&\lesssim \left(\frac{1}{v(D(a_i,2r))}\int_{D(a_i,2r)}
					|f-h_i|^q vdA\right)^\frac{1}{q}\\
					&+\left(\frac{1}{v(D(a_j,2r))}\int_{D({a_j,2r})}|f-h_j|^qvdA \right)^\frac{1}{q}\\
					&\lesssim G_{q,2r}(f)(a_i)+G_{q,2r}(f)(a_j)+\frac{1}{i}+\frac{1}{j}\\
					&\lesssim G_{q,2r}(f)(w)+{C}.  
				\end{aligned}
			\end{align}

			Let \(I_z\) be the set of integers \(i\) such that \(z\in{D(a_i,r)}\) and \(j\) be an element in \(I_z\).  Set \[f_1=\sum\limits_i h_i\phi_i=h_j+\sum_i (h_i-h_j)\phi_i,\]
			then \[\overline{\partial}f_1=\sum_{i}(h_i-h_j)\overline{\partial}\phi_i.\] 
			There are at most \(N\) of the integers \(i\) such that \(\text{Supp}\,\phi_i\ni z\), combing with \eqref{eqhihj} yields 
			\begin{equation}\label{eqf1}
				|\overline{\partial}f_1(z)|\leq\frac{CN}{1-|z|}G_{q,2r}(f)(z),
			\end{equation}
			hence,
			\[[v]^{-\frac{1}{p}}[\eta]^\frac{1}{q}\lesssim[\eta]^\frac{1}{q}[v]^{-\frac{1}{p}}G_{q,2r}(f)<\infty.\]

			Let \(f_2=f-f_1\), 
			\begin{align*}
				\left(\frac{1}{|D(z,r)|}\int_{D(z,r)}|f_2|^qdA\right)^\frac{1}{q}&=\left(\sum_i\frac{1}{|D(z,r)|}\int_{D(z,r)}|f-h_i|^q|\phi_i|^qdA\right)^\frac{1}{q}\\
				&\leq\sum_i\left(\frac{1}{|D(z,r)|}\int_{D(z,r)\bigcap D(a_i,r)}|f-h_i|^qdA\right)^\frac{1}{q}.\\
			\end{align*}
			Since \(D(z,r)\bigcap D(a_i,r)\neq\emptyset\) and the fact that \(D(z,r)\bigcap D(a_i,r)\) can be covered by at most \(N\) integers \(i\), then by Lemma \ref{weights} yields 
			\begin{equation}\label{Mrf2}
				M_r(|f_2|^q)^\frac{1}{q}\lesssim G_{q,r}(f)(z).\\
			\end{equation}
			It implies that
			\[[\eta]^\frac{1}{q}[v]^{-\frac{1}{p}}M_r(|f_2|)^\frac{1}{q}<\infty.\]

			\((\text{iii})\implies (\text{i})\).
			Write \(f=f_1+f_2\), \(H^\om_{f}(g)=H^\om_{f_1}(g)+H^\om_{f_2}(g)\), \(d\nu=(|1-|z|)|\overline{\partial}f_1|\eta dA\).
			Then by Lemma \ref{CMpq}, Lemma \ref{dbar} and \cite[Theorem 3.2]{Hu2019},
			\begin{align*}
				\|H^\om_{f_1}(g)\|_{L_\eta^q}&=\|u-P_\om(u)\|_{L^q_\eta}\\
				&\leq(1+\|P_\om\|)\|u\|_{L_\eta^q}\\
				&\leq C\|(1-|z|)g\overline{\partial}f_1\|_{L_\eta^q}\\
				&\leq\frac{\int_{D(z,r)}d\nu}{v(D(z,r))^\frac{q}{p}}\|g\|_{A_v^p}\\
				&\asymp[\eta][v]^{-\frac{q}{p}}G^q_{q,r}(f)\|g\|_{A_v^p}\\
				&\lesssim \|g\|_{A_v^p}.
			\end{align*}
			For \(|f_2|^q\eta dA\), by Lemma \ref{weights}, 
			\begin{align*}
				\| H_{f_2}^\om(g)\|_{L_\eta^q}&=\|(Id-P_\om)(f_2g)\|_{L_\eta^q}\\
				&\leq \|f_2g\|_{L_\eta^q}\\
				&\lesssim \frac{\int_{D(z,r)}|f_2|^q\eta dA}{v(D(z,r))^\frac{q}{p}}\|g\|_{A_v^p}\\
				&\asymp [v][\eta]^{-\frac{1}{p}}M_r(|f_2|^q)\\
				&\lesssim \|g\|_{A_v^p}.
			\end{align*}
			The proof is completed.
		\end{proof}
		
		The characterization for the compactness of Hankel operator is obtained as follows. 
		\begin{theorem}\label{mainth2}
			Let \(\om,v,\eta\in\cR\), \(1<p\leq q<\infty\),  \(A_p(\om,\eta)<\infty\), \(f\in L_v^1\). Then the following statements are equivalent. 
			\begin{enumerate}
				\item [(i)] \(H_f^\om: A_v^p\to L_\eta^q\) is compact;
				\item [(ii)] For some (or any) \(r>0\), \(\lim\limits_{|z|\to 1^-}[\eta]^\frac{1}{q}[v]^{-\frac{1}{p}}G_{q,r}(f)(z)=0\);
				\item [(iii)] \(f\) admits a decomposition \(f=f_1+f_2\), where \(f_1\in C^1(\bbD)\) and satisfies 
				\[\lim_{|z|\to 1^-|}[\eta]^\frac{1}{q}[v]^{-\frac{1}{p}}(1-|z|)|\overline{\partial}f_1|=0,\]
				and \(f_2\) satisfies 
				\[\lim_{|z|\to 1^-}[\eta]^\frac{1}{q}[v]^{-\frac{1}{p}}M_r(|f_2|^q)^\frac{1}{q}=0.\]
			\end{enumerate}
		\end{theorem}
		\begin{proof}
			The proof of Theorem \ref{mainth2} is similar to that of the Theorem \ref{mainth1}, so we only provide a general outline.  Suppose that \(H_f^\om\) is compact. Since \(b_z^\om\to 0\) weakly, then \eqref{eqHf} yields 
			\[\lim_{|z|\to 1^-}[\eta]^\frac{1}{q}[v]^{-\frac{1}{p}}G_{q,r}(f)(z)=0.\]
			This is conclusion (ii).

			Suppose (ii) holds. Then by \eqref{Mrf2} and \eqref{eqf1} deduce (iii) holds. 
			
			Finally, we suppose that (iii) holds. Let \(\{g_n\}\) converge to zero uniformly on any compact subset of \(\bbD\), then we need only show that \(\|H_f^\om(g_n)\|_{L_\eta^q}\to 0\) as \(n\to\infty\). From the proof of Theorem \ref{mainth1}, \((|1-|z|)|\overline{\partial}f_1|\eta dA\) and \(|f_2|^q\eta dA\) are vanishing \(q\)-Carleson measure, then by Lemma \ref{weights} and Lemma \ref{CMpq} we get (i). The proof is completed.
		\end{proof}
		
		\section{The case of \(1<q<p<\infty\)}
		The characterization for the case of \(1<q<p<\infty\) to Hankel operators is obtained  by using Khintchine's inequality in this section.
		\begin{theorem}\label{mainth3}
			Let \(1<q<p<\infty\), \(\om,v,\eta\in\cR\),  \(A_p(\om,\eta)<\infty\), \(W=\eta^\frac{p}{p-q}v^{-\frac{q}{p-q}}\). Then for \(f\in L_v^1\), the following statements are equivalent. 
			\begin{enumerate}
				\item [(i)] \(H_f^\om:A_v^p\to L_\eta^q\) is compact;
				\item [(ii)] \(H_f^\om:A_v^p\to L_\eta^q\) is bounded;
				\item [(iii)] For some (or any) \(r\), \(G_{q,r}(f)(z)\in L_W^\frac{pq}{p-q}\);
				\item [(iv)] \(f\) admits a decomposition \(f=f_1+f_2\), where \(f_1\in C^1(\bbD)\)  satisfies 
				\((1-|z|)|\overline{\partial}f_1|\in L_W^\frac{pq}{p-q}\), and \(f_2\) satisfies
				\(M_r(|f_2|^q)^\frac{1}{q}\in L_W^\frac{pq}{p-q}.\)
			\end{enumerate}
		\end{theorem}
		\begin{proof}
			We prove the theorem in its natural order.
			
			\((\text{i})\implies (\text{ii})\) is trivial.
			
			\((\text{ii})\implies (\text{iii})\).
			Suppose that  \(\{\lambda_j\}_{j=1}^\infty\in l^p\), \(\{z_j\}\) is an \(\frac{r}{4}\)-lattice and \(b_{v,z_j}=\frac{B_{z_j}^\om(z)}{\|B_{z_j}^\om\|_{A_v^p}}\),  set \[F(z)=\sum_j\lambda_j\frac{B_{z_j}^\om(z)}{\|B_{z_j}^\om\|_{A_v^p}}=\sum_j\lambda_jb_{v,z_j}.\] 
			Let \(r_j(t)\) be a Radermacher random variables on \(L^1((0,1),dt)\), \(R_F=\sum_{j=1}^\infty r_j(t)\lambda_jb_{v,z_j}\) is the Radermacher randomization of \(F\). Then Fubini's theorem, Khintchine's inequality and Lemma \ref{normestimates} yields 
			\begin{align*}
				&\int_0^1\|H_f^\om(R_F)\|_{L_\eta^q}^qdt\\
				&=\int_0^1\left\|\sum_{j}\lambda_jr_j(t)H^\om_f(b_{v,z_j})\right\|_{A_v^p}^qdt\\
				&\asymp\int_\bbD\left(\sum_j|\lambda_j|^2|H^\om_f(b_{v,z_j})^2|\right)^\frac{q}{2}\eta(z)dA(z)\\
				&=\sum_{k}\int_{D(z_k,r)}\left(\sum_j|\lambda_j|^2|H^\om_f(b_{v,z_j})^2|\right)^\frac{q}{2}\eta(z)dA(z)\\
				&\gtrsim\sum_k|\lambda_k|^q|b_{v,z_k}(z_k)|^q\int_{D(z_k,r)}\left|fb_{v,z_j}-\frac{P_w(b_{v,z_j})}{b_{v,z_j}}\right|^q\eta(z)dA(z)\\
				&\gtrsim\sum_k |\lambda_k|^qG_{q,r}^q(f)(z_k)\eta(D(z_k,r))v(D(z_k,r))^{-\frac{q}{p}}.
			\end{align*} 
			Since \[\left\|\sum_j\lambda_jb_{v,z_j}\right\|_{A_v^p}\lesssim \|\lambda_j\|_{l^p}=\||\lambda_j|^q\|_{l^\frac{p}{q}}^\frac{1}{q},\]
			and 
			\begin{align*}
				\left\|H_f^\om\left(\sum_j\lambda_jr_j(t)b_{v,z_j}\right)\right\|_{L_\eta^q}\lesssim \left\|\sum_j\lambda_jr_j(t)b_{v,z_j}\right\|_{A_v^p}\asymp\||\lambda_j|^q\|_{l^\frac{p}{q}}^\frac{1}{q},
			\end{align*}
			then 
			\begin{align*}
				\sum_{k}|\lambda_k|^q\eta(D(z_k,r))v(D(z_k,r))^{-\frac{q}{p}}G_{q,r}^q(f)(z_k)\lesssim \|H_f^\om\|\|\lambda_j\|_{l^\frac{p}{q}}^\frac{1}{q}.
			\end{align*}
			Notice that \[\frac{1}{\frac{p}{q}}+\frac{1}{\frac{p}{p-q}}=1,\] an argument of duality yields 
			\begin{align*}
				\sum_kW(D(z_k,r))G_{q,r}^\frac{pq}{p-q}(f)(z_k)\lesssim \|H_f^\om\|.
			\end{align*}
			For the point \(z\) and \(w\) satisfy \(\beta(z,w)<\frac{r}{2}\), 
			\begin{equation}\label{eqasymp}
				G_{q,r}(f)(w)\asymp G_{q,\frac{r}{2}}f(z).
			\end{equation}
			Hence,
			\begin{align*}
				\infty> \|H_f^\om\|&\gtrsim\sum_kW(D(z_k,\frac{r}{2}))G_{q,r}^\frac{pq}{p-q}(f)(z_k)\\
				&\gtrsim\sum_k\int_{(D(z_k,\frac{r}{2}))}G_{q,r}^\frac{pq}{p-q}(z)W(z)dA(z)\\
				&\gtrsim\int_\bbD G_{q,r}^\frac{pq}{p-q}(f)(z)W(z)dA(z).
			\end{align*}

			\((\text{iii})\implies (\text{iv})\).
			Let \(\phi_j\) be a smooth partition of \(\bbD\) subordinate to \(D(a_j,r)\), set 
			\(f_2=f-f_1\) and \(f_1=\sum\limits_jh_j\phi_j\in C^\infty(\bbD)\) as in the proof of Theorem \ref{mainth1}.  By \eqref{eqf1},
			\[|\overline{\partial}f_1(z)|\leq\frac{CN}{1-|z|}G_{q,\frac{r}{2}}(f)(z).\]
			Hence, 
			\begin{align*}
				\left((1-|z|)|\overline{\partial}f_1(z)|\right)^\frac{pq}{p-q}&\lesssim\frac{1}{W(D(a_j,\frac{r}{2}))}\int_{D(a_j,r)}G_{q,r}^\frac{pq}{p-q}(z)(f)(z)W(z)dA(z) \\
				&\lesssim \frac{1}{W(D(z,r))}\int_{D(z,2r)}G_{q,r}^\frac{pq}{p-q}(f)(w)W(w)dA(w).
			\end{align*}
			By Fubini's theorem, 
			\begin{align*}
				&\int_\bbD \left((1-|z|)|\overline{\partial}f_1(z)|\right)^\frac{pq}{p-q}W(z)dA(z)\\
				&\lesssim\int_\bbD\frac{1}{W(D(z,r))}\int_{D(z,2r)}G_{q,r}^\frac{pq}{p-q}(f)(w)W(w)dA(w)W(z)dA(z)\\
				&=\int_\bbD\frac{1}{W(D(z,r))}\chi_{D(w,2r)}(z)W(z)dA(z)\int_{D(z,2r)}G_{q,r}^\frac{pq}{p-q}(f)(w)W(w)dA(w)\\
				&\asymp\int_\bbD G_{q,r}^\frac{pq}{p-q}W(u)dA(u).
			\end{align*}

			Similar to the proof of \eqref{Mrf2} and the fact \eqref{eqasymp} yields
			\begin{align*}
				M_r(|f_2|^q)^\frac{1}{q}(z)\lesssim G_{q,2r}(f)(z).
			\end{align*}
			Hence,
			\begin{equation*}
				\int_\bbD \left(M_r(|f_2|^q)^\frac{1}{q}(z)\right)^\frac{pq}{p-q}W(z)dA(z)\lesssim\int_\bbD G_{q,2r}^\frac{pq}{p-q}(f)(z)W(z)dA(z)<\infty. 
			\end{equation*}
			As in the proof of Theorem \ref{mainth1}, write \(d\nu=(1-|z|)|\overline{\partial}f_1|^q\eta dA,\)  then by Lemma \ref{CMqp} and H\"older's inequality yields 
			\begin{align*}
				&\|H_{f_1}^\om(g)\|_{L_\eta^q}\\
				&\leq(1+\|P_\om\|)\|(u)\|_{L_\eta^q}\\
				&\lesssim\|(1-|z|)g\overline{\partial}f_1\|_{L_\eta^q}\\
				&\lesssim \int_\bbD\left(\frac{\int_{D(z,r)}(1-|\zeta|)^q|\overline{\partial}f_1|^q\eta(\zeta)dA(\zeta)}{v(D(z,r))}\right)^\frac{p}{p-q}v(\zeta)dA(\zeta)\|g\|_{A_v^p}\\
				&\lesssim\int_\bbD\left(\int_{D(\zeta,r)}((1-|\zeta|)|\overline{\partial}f_1|)^\frac{pq}{p-q}\eta(\zeta)dA(\zeta)\right)\\
				&\times\left(\int_{D(z,r)}\left(v(D(z,r))\right)^{-\frac{p}{q}}\eta(\zeta)dA(\zeta)\right)^\frac{q}{p-q}v(z)dA(z)\|g\|_{A_v^p}\\
				&=\int_\bbD\left(\int_{D(\zeta,r)}((1-|\zeta|)|\overline{\partial}f_1|)^\frac{pq}{p-q}W(\zeta)\eta(\zeta)^{-\frac{q}{p-q}}v(\zeta)^\frac{p}{p-q}\frac{1}{v(\zeta)}dA(\zeta)\right)\\
				&\times\eta(D(z,r))^{\frac{q}{p-q}}v(D(z,r))^{-\frac{p}{p-q}}\eta(z)dA(z)\|g\|_{A_v^p}\\
				&\asymp\int_\bbD \left((1-|z|)|\overline{\partial}f_1|\right)^\frac{pq}{p-q}W(z)dA(z)\|g\|_{A_v^p}.
			\end{align*}
			Similarly, write \(d\mu=|f_2|^q\eta dA,\)
			then 
			\begin{align}\label{eq11}
				\begin{aligned}
					\|H_{f_2}^\om(g)\|_{L_\eta^q}&\leq\|f_2g\|_{L_\eta^q}\\
					&\lesssim\int_\bbD M_r(|f_2|^q)^\frac{1}{q}(z)W(z)dA(z)\|g\|_{A_v^p}.
				\end{aligned}
			\end{align}
			
			\((\text{iv})\implies (\text{i})\).
			Next we show that \(H_{f_1}^\om\) is compact. Suppose that \(\{g_m\}\) converges to zero uniformly on any subsets of \(\bbD\), we need to show that \(H_{f_1}^\om(g_m)\) converges to zero strongly. Since \(A_v^p\) is Banach space, then there exists \(h_{m,n}\) such that 
			\[\|g_m-h_{m,n}\|_{A_v^p}<\frac{1}{m}.\]
			
			Let 
			\begin{equation*}
				u_m(z)=\sum_jB_{z_j}(z)\int_\bbD\frac{\phi_j(\zeta)}{(\zeta-z)B_{z_j}(z)}h_{m,n}(\zeta)\overline{\partial}f_1(\zeta)dA(\zeta),
			\end{equation*}
			then \(\overline{\partial}u_m=h_{m,n}\overline{\partial}f_1\)
			and  $ \|u_m\|_{L_\eta^q}\lesssim \|h_{m,n}\overline{\partial}f_1\|_{L^q((1-|z|)^q\eta)}$.
			It follows from Lemma \ref{CMqp} that \(Id:A_v^p\hookrightarrow L^q(\nu)\) is compact. Hence 
			\begin{align*}
				\|u_m\|_{L_\eta^q}&\lesssim \|h_{m,n}\overline{\partial}f_1\|_{L^q((1-|z|)^q\eta)}\to 0, \,\,m\to\infty.
			\end{align*}
			Since 
			\(H_{f_1}^\om(g_{m,n})=u_m-P_\om(u_m),\)
			then 
			\begin{align*}
				\|H_{f_1}^\om(h_{m,n})\|_{L_\eta^q}\leq (1+\|P_\om\|)\|u\|_{L_\eta^q}\to 0,\,\,m\to\infty.
			\end{align*}
			Since 
			\begin{align*}
				\|H_{f_1}^\om(g_m-g_{m,n})\|_{L_\eta^q}\leq\|H_{f_1}^\om\|\|g_m-g_{m,n}\|_{A_v^p}\to 0,\,\, m\to\infty,
			\end{align*}
			then 
			\begin{align*}
				\|H_{f_1}^\om(g_m)\|_{L_\eta^q} \leq\|H_{f_1}^\om(g_m-g_{m,n})\|_{L_\eta^q}+\|H_{f_1}^\om(g_{m,n})\|_{L_\eta^q}\to 0,\,\,m\to\infty.
			\end{align*}
			Similarly, write \(d\mu=|f_2|^q\eta dA\), then by Lemma \ref{CMqp} and \eqref{eq11},
			\begin{align*}
				\|H_{f_2}(g_m)\|_{L_\eta^q}&\leq \|f_2g\|_{L_\eta^q}\\
				&\lesssim\|Id\|_{A_v^p\rightarrow L^q(d\mu)}\|g_m\|_{A_v^p}\to 0,\,\,m\to\infty.
			\end{align*}
			Therefore, \(H_f^\om\) is compact. The proof is completed.
		\end{proof}
		
		\section{Acknowledgement}
		This research work is supported by the National Natural Science Foundation of China~(Grant No.~12261023,~11861023.)


\begin{thebibliography}{99}
			\bibitem{Aleman2017} A. Aleman, S. Pott,   M. Reguera, Sarason Conjecture on the Bergman Space,  Int. Math. Res. Not.,  14 (2017),  4320--4349.
			\bibitem{Axler1986} S. Axler, The Bergman space, the Bloch space and commutators of multiplication operators, Duke. Math. J., 53 (1986), 315-332.
			\bibitem{Arazy1988} J. Arazy, S. Fisher, J. Peetre, Hankel operators on weighted Bergman spaces, Amer. J. Math., 110 (1988), 989-1054.
			\bibitem{Bekolle1990} D. Bekoll\'e, C. Berger, L. Coburn, K. Zhu, \(BMO\) in the Bergman metric on bounded symmetric domains, J. Funct. Anal., 93 (1990), 310-350.
			\bibitem{Duan2022} Y. Duan, K. Guo, S. Wang, Z. Wang, Toeplitz operators on weighted Bergman spaces induced by a class of radial weights, J. Geom. Anal., 32 (2022), Article number 39, 29pp.
			\bibitem{Duan2023} Y. Duan, J. R\"atty\"a, S. Wang, F. Wu, Two weight inequality for Hankel form on weighted Bergman spaces induced by doubling weights, Adv. Math., 431 (2023), Paper ID: 109249.
			\bibitem{Hu2018} Z. Hu, E. Wang, Hankel operators between Fock spaces, Integr. Equ. Oper. Theory, 90 (2018), Article number 37, 20pp.
			\bibitem{Hu2019} Z. Hu, J. Lu, Hankel operators on Bergman spaces with regular weights, J. Geom. Anal., 29 (2019), 3494-3519.
			\bibitem{Hu2023} Z. Hu, J. Virtanen, IDA and Hankel operators on Fock spaces, Anal. and PDE, 16 (2023), 2041-2077. 
			\bibitem{Hytonen2012} T. P. Hyt\"{o}nen, The sharp weighted bound for general Calder\'{o}n-Zygmund operators, Ann. Math.,  175  (2012), no. 3, 1473--1506.
			\bibitem{Jason1988} S. Jason, Hankel operators between Bergman spaces, Ark. Mat., 26 (1988), 205-219.
			\bibitem{Lacey2014} M.  Lacey, E. Sawyer, C. Shen,  I. Uriarte-Tuero, Two-weight inequality for the Hilbert transform: A real variable characterization, I,  Duke Math. J., 163 (2014),  2795--2820.
			\bibitem{Luecking1992} D.  Luecking, Characterizations of certain classes of Hankel operators on the Bergman spaces of the unit disc, J. Funct. Anal., 110 (1992), 247-271.
			\bibitem{martin2021} F. Mart\'in Reyes, P. Ortega, J. Pel\'aez, J. R\"atty\"a, One weight inequality for Bergman projection and Calder\'on operators induced by radial weights, arxiv: 2105.08029v1 (2021).
			\bibitem{Pau2013} J. Pau, Hankel operators on standard Bergman spaces, Complex Anal. Oper. Theory, 7 (2013), 1239-1256.
			\bibitem{Pau2016} J. Pau, R. Zhao, K. Zhu, Weighted \(BMO\) and Hankel operators between Bergman spaces, Indiana Univ. Math. J., 65 (2016), 1639-1673. 
			\bibitem{Pelaez2013}
			J.  Pel\'aez, J. R\"atty\"a, Generalized Hilbert operators on weighted Bergman spaces, Adv. Math., 240 (2013), 227-267.
			\bibitem{Pelaez2014}  J.  Pel\'aez, R\"atty\"a,  Weighted Bergman spaces induced by rapidly increasing weights, Mem. Amer. Math. Soc., 227 (1066), 124 (2014).
			\bibitem{Pelaez2015}
			J.  Pel\'aez, J. R\"atty\"a, Embedding theorems for Bergman spaces via harmonic analysis, Math Ann., 362 (2015), 205-239.
			\bibitem{Pelaez2016summer} J.  Ple\'aez, Small weighted Bergman spaces, Proceedings of the Summer school in Complex and Harmonic Analysis, and Related Topics, University of Eastern Finland, Faculty of Science and Forestry, Jonesuu, 2016, 29-98.
			\bibitem{Pelaezprojection2016}
			J.  Pel\'aez, J. R\"atty\"a, Two inequality for Bergman projection, J. Math. Pures Appl., 105 (2016), 102-130.
			\bibitem{Pelaez2018}
			J.  Pel\'aez, J. R\"atty\"a, K. Sierra, Berezin transform and Toeplitz operators on Bergman spaces induced by regular weights,  {J. Geom. Anal.} 28 (2018), 656-687.
			\bibitem{Pelaez2021}
			J.  Pel\'aez, J. R\"atty\"a, Bergman projection induced by radial weight, Adv. Math., 70 (2021), Paper ID: 107950.
			\bibitem{Pelaez2020} J.  Pel\'aez, A. Per\"al\"a, J. R\"atty\"a, Hankel operators induced by radial Bekoll\'e-Bonami weights on Bergman spaces, Math. Z., 296 (2020), 211-238.
			\bibitem{Stroethoff1990} K. Stroethoff, Compact Hankel operators on the Bergman space, Illinois J. Math., 34 (1990), 159-174.
			\bibitem{Korhonen2021} T. Korhonen, J. A. Pel\'aez, J. R\"atty\"a, Radial two weight inequality for Maximal Bergman Projection induced by a regular weight, Potential Anal., 54 (2021), 561-574.
			\bibitem{Zhang2025} P. Zhang, M. Li, J. Long, On two weights inequality induced by radial weights on the unit ball of \(\bbC^{n}\), to appear in Southeast Asian Bull. Math.
			\bibitem{Zhu2007}
			K. Zhu,  {Operator Theory in Function Spaces}, second ed., Math. Surveys Monogr., vol. 138, American Mathematical Society, Providence, RI, 2007.
			
		\end{thebibliography}
	\end{document}